\documentclass[12pt,a4paper,twoside]{article}
\usepackage[english]{babel}
\usepackage{amssymb}
\usepackage{inputenc}
\newcommand{\ba}{\begin{array}}
\newcommand{\ea}{\end{array}}
\usepackage{amssymb}
\begin{document}

\author{S. Albeverio $^{1},$ Sh. A. Ayupov $^{2,  *},$  K. K.
Kudaybergenov  $^3,$  T. S. Kalandarov  $^4$}
\title{\bf Complete description of derivations  on $\tau$-compact
 operators for type I  von Neumann algebras}

\maketitle

\begin{abstract}

Given a type I von Neumann algebra  $M$ with a faithful normal
semi-finite trace $\tau,$ let $S_0(M, \tau)$ be the algebra of all
$\tau$-compact operators affiliated with $M.$ We give a complete
description of all derivations on the algebra $S_0(M, \tau).$ In
particular, we prove that if $M$ is of type I$_\infty$ then every
derivation on $S_0(M, \tau)$ is spatial.

\end{abstract}

\medskip
$^1$ Institut f\"{u}r Angewandte Mathematik, Universit\"{a}t Bonn,
Wegelerstr. 6, D-53115 Bonn (Germany); SFB 611, BiBoS; CERFIM
(Locarno); Acc. Arch. (USI), \emph{albeverio@uni-bonn.de}

$^2$ Institute of Mathematics and information  technologies,
Uzbekistan Academy of Sciences, F. Hodjaev str. 29, 100125,
Tashkent (Uzbekistan), e-mail: \emph{sh\_ayupov@mail.ru}

 $^{3}$ Karakalpak state university, Ch. Abdirov str. 1, 742012, Nukus (Uzbekistan),
e-mail: \emph{karim2006@mail.ru}

 $^{4}$ Institute of
Mathematics and information  technologies, Uzbekistan Academy of
Sciences, F. Hodjaev str. 29, 100125, Tashkent (Uzbekistan),
e-mail: \emph{turaboy\_kts@rambler.ru}



* Corresponding author
\newpage

\begin{center}
{\bf 1. Introduction}
\end{center}

Derivations on unbounded operator algebras, in particular on
various algebras of measurable operators affiliated with von
Neumann algebras, appear to be a very attractive special case of
the general theory of unbounded derivations on operator algebras.
The present paper continues the series of papers of the authors
\cite{Alb1}-\cite{Alb4} devoted to the study and a description of derivation on the
algebra $LS(M)$ of locally measurable operators with respect to a
von Neumann algebra $M$ and on various subalgebras of $LS(M).$

Let  $A$ be an algebra over the complex number. A linear operator
$D:A\rightarrow A$ is called a derivation if it satisfies the
identity  $D(xy)=D(x)y+xD(y)$ for all  $x, y\in A$ (Leibniz rule).
Each element  $a\in A$ defines a derivation  $D_a$ on $A$ given as
$D_a(x)=ax-xa,\,x\in A.$ Such derivations $D_a$ are said to be
\emph{inner derivations}. If the element  $a$ generating the
derivation   $D_a$ on $A,$ belongs to a larger algebra  $B,$
containing  $A$ (as a proper ideal as usual) then $D_a$ is called
a  \emph{spatial derivation}.

In the particular case where $A$ is commutative, inner derivations
are identically zero, i.e. trivial. One of the main problems in
the theory of derivations is automatic innerness or spatialness of
derivations and the existence of non inner derivations (in
particular, non trivial derivations on commutative algebras).

On this way  A.~F.~Ber, F.~A.~Sukochev, V.~I.~Chilin~\cite{Ber}
obtained necessary and sufficient conditions for the existence of
non trivial derivations on commutative regular algebras. In
particular they have proved that the algebra  $L^{0}(0, 1)$ of all
(classes of equivalence of) complex measurable functions on  the interval $(0,
1)$   admits non trivial derivations. Independently
A.~G.~Kusraev~\cite{Kus1} by means of Boolean-valued analysis has
established necessary and sufficient conditions for the existence
of non trivial derivations and automorphisms on extended
$f$-algebras. In particular he has also proved the existence of
non trivial derivations and automorphisms on  $L^{0}(0, 1).$ Later
the authors initiated the study of these problems in the non
commutative case  \cite{Alb1}-\cite{Alb4}, by considering derivations
on the algebra $LS(M)$ of all locally measurable operators with
respect to a semi-finite von Neumann algebra, and on various
subalgebras of  $LS(M).$ Thus in  \cite{Alb1}, \cite{Alb2} we
obtained  necessary and sufficient conditions for derivations to
be inner on the algebra $LS(M)$ and its subalgebra $S_0(M, \tau)$
of  $\tau$-compact operators with respect to  $M,$ where  $M$ is a
von Neumann algebra of type I with a faithful normal semi-finite
trace  $\tau.$  Moreover in  \cite{Alb4} we gave a complete
descriptions of all derivations on the algebra $LS(M),$ for type I
von Neumann algebra $M.$ In \cite{Alb3} we have proved spatialness
of all derivations on the non commutative Arens algebra
$L^{\omega}(M, \tau)$ associated with an arbitrary von Neumann
algebra  $M$ and a faithful normal semi-finite trace  $\tau.$

In the present paper which extends  \cite{Alb4} we present a
complete descriptions of derivations on the algebra  $S_0(M,
\tau)$ of  $\tau$-compact operators affiliated with a type I von
Neumann algebra $M$ and a faithful normal semi-finite trace
$\tau.$

\begin{center}
{\bf 2. Preliminaries}
\end{center}

Let  $H$ be a complex Hilbert space and let  $B(H)$ be the algebra
of all bounded linear operators on   $H.$ Consider a von Neumann
algebra $M$  in $B(H)$ with a faithful normal semi-finite trace
$\tau.$ Denote by  $P(M)$ the lattice of projections in $M.$

A linear subspace  $\mathcal{D}$ in  $H$ is said to be
\emph{affiliated} with  $M$ (denoted as  $\mathcal{D}\eta M$), if
$u(\mathcal{D})\subset \mathcal{D}$ for every unitary  $u$ from
the commutant
$$M'=\{y\in B(H):xy=yx, \,\forall x\in M\}$$ of the von Neumann algebra $M.$

A linear operator  $x$ on  $H$ with the domain  $\mathcal{D}(x)$
is said to be \emph{affiliated} with  $M$ (denoted as  $x\eta M$) if
$\mathcal{D}(x)\eta M$ and $u(x(\xi))=x(u(\xi))$
 all  $\xi\in
\mathcal{D}(x).$

A linear subspace $\mathcal{D}$ in $H$ is said to be \emph{strongly
dense} in  $H$ with respect to the von Neumann algebra  $M,$ if

1) $\mathcal{D}\eta M;$

2) there exists a sequence of projections
$\{p_n\}_{n=1}^{\infty}$ in $P(M)$  such that
$p_n\uparrow\textbf{1},$ $p_n(H)\subset \mathcal{D}$ and
$p^{\perp}_n=\textbf{1}-p_n$ is finite in  $M$ for all
$n\in\mathbb{N},$ where $\textbf{1}$ is the identity in $M.$

A closed linear operator  $x$ acting in the Hilbert space $H$ is said to be
\emph{measurable} with respect to the von Neumann algebra  $M,$ if
 $x\eta M$ and $\mathcal{D}(x)$ is strongly dense in  $H.$ Denote by
 $S(M)$ the set of all measurable operators affiliated with
 $M.$

A closed linear operator $x$ in  $H$  is said to
be \emph{locally measurable} with respect to the von Neumann
algebra $M,$ if $x\eta M$ and there exists a sequence
$\{z_n\}_{n=1}^{\infty}$ of central projections in $M$ such that
$z_n\uparrow\textbf{1}$ and $z_nx \in S(M)$ for all
$n\in\mathbb{N}.$

It is well-known \cite{Mur1} that the set $LS(M)$ of all locally
measurable operators with respect to $M$ is a unital *-algebra when equipped with
 the algebraic operations of strong addition and multiplication and taking the adjoint of an operator.

 Let   $\tau$ be a faithful normal semi-finite trace on $M.$ Recall that a closed linear operator
  $x$ is said to be  $\tau$\emph{-measurable} with respect to the von Neumann algebra
   $M,$ if  $x\eta M$ and   $\mathcal{D}(x)$ is
  $\tau$-dense in  $H,$ i.e. $\mathcal{D}(x)\eta M$ and given   $\varepsilon>0$
  there exists a projection   $p\in M$ such that   $p(H)\subset\mathcal{D}(x)$ and $\tau(p^{\perp})<\varepsilon.$
   The set $S(M,\tau)$ of all   $\tau$-measurable operators with respect to  $M$
    is a solid  *-algebra in $S(M)$  (see \cite{Mur1}).

An element  $x$ of the algebra  $S(M, \tau)$ is said to be  $\tau$-\emph{compact}, if given any
 $\varepsilon >0$ there exists a projection $p\in P(M)$
such that  $\tau(p^{\perp})<\infty,\,xp\in M$ and
$\|xp\|_M<\varepsilon.$ The set $S_0(M, \tau)$ of all $\tau$-compact operators is a  *-ideal in the algebra
$S(M, \tau)$
(see \cite{Mur1}).

It is well-known \cite{Seg} that every commutative von Neumann algebra
 $M$
is *-isomorphic to the algebra   $L^{\infty}(\Omega)=L^{\infty}(\Omega, \Sigma,
\mu)$ of all (classes of equivalence of) complex essentially bounded measurable functions on a measure space
$(\Omega, \Sigma, \mu)$ and in this case  $LS(M)=S(M)\cong
L^{0}(\Omega),$   where $L^{0}(\Omega)=L^{0}(\Omega, \Sigma,
\mu)$ the algebra of all (classes of equivalence of) complex measurable functions on
$(\Omega, \Sigma, \mu).$

Now let us recall some notions and results from the theory of
Hilbert -- Kaplansky modules (for details we refer to \cite{Kus2}).

Let    $H$ be a Hilbert space and denote by  $L^{0}(\Omega, H)$
 the space of all equivalence classes of measurable maps from
 $\Omega$  into $H.$
Equipped with the  $L^{0}(\Omega)$-valued inner product
   $$\langle x, y \rangle=\langle x(\omega), y(\omega) \rangle_{H},$$ where  $\langle \cdot, \cdot\rangle_{H}$
   in the inner product in  $H,$
      $L^{0}(\Omega, H)$ becomes a Hilbert~-- Kaplansky module over $L^{0}(\Omega).$
       The space   $$L^{\infty}(\Omega, H)=\{x\in L^{0}(\Omega, H): \langle x,x\rangle\in L^{\infty}(\Omega)\}$$
      is a Hilbert~-- Kaplansky module over  $L^{\infty}(\Omega).$ Denote by   $B(L^{0}(\Omega, H))$ the
      algebra of all
              $L^{0}(\Omega)$-bounded
       $L^{0}(\Omega)$-linear operators on  $L^{0}(\Omega, H)$ and denote by
        $B(L^{\infty}(\Omega, H))$ the algebra of all
        $L^{\infty}(\Omega)$-bounded $L^{\infty}(\Omega)$-linear operators on
        $L^{\infty}(\Omega, H).$

Now consider a von Neumann algebra  $M$ which is homogeneous of type
I$_{\alpha}$ with the center  $L^{\infty}(\Omega),$ where $\alpha$ is a cardinal number.
Then    $M$ is *-isomorpic to the algebra
  $B(L^{\infty}(\Omega, H)),$ where $\dim H=\alpha,$ while the algebra  $LS(M)$
  is *-isomorpic to $B(L^{0}(\Omega, H))$ (see for details \cite{Alb1}).

It is known  \cite{Tak} that given a type I von Neumann algebra
$M$ there exists a unique (cardinal-indexed) family  of central
orthogonal projections  $(q_{\alpha})_{\alpha\in J}$ in
$\mathcal{P}(M)$ with $\sum\limits_{\alpha\in
J}q_{\alpha}=\textbf{1}$ such that $q_{\alpha}M$ is a homogeneous
type I$_\alpha$ von Neumann algebra, i.e.
 $q_\alpha
M\cong B(L^{\infty}(\Omega_{\alpha}, H_{\alpha}))$ with $\dim
H_\alpha=\alpha$  and
$$M\cong\bigoplus\limits_{\alpha\in J}B(L^{\infty}(\Omega_{\alpha}, H_{\alpha}))$$
($\cong$ denoting isomorphism).
The direct product
$$\prod\limits_{\alpha\in J}L^{0}(\Omega_{\alpha}, H_{\alpha})$$
equipped with the coordinate-wise algebraic operations and inner product forms
a Hilbert~--- Kaplansky module over $L^{0}(\Omega).$

In  \cite{Alb1} we have proved that for
$M\cong\bigoplus\limits_{\alpha\in J}B(L^{\infty}(\Omega_{\alpha},
H_{\alpha}))$ the algebra   $LS(M)$ is  *-isomorpic with   $B(\prod\limits_{\alpha\in J}L^{0}(\Omega_{\alpha},
H_{\alpha})).$ Therefore there exists a map $||\cdot||:LS(M)\rightarrow L^{0}(\Omega)$ such that for all  $x,y\in
LS(M), \lambda\in L^{0}(\Omega)$ one has
$$||x||\geq0, ||x||=0 \Leftrightarrow x=0;$$
$$||\lambda x||=|\lambda|||x||;$$
$$||x+y||\leq||x||+||y||;$$
$$||xx^{*}||=||x||^{2};$$
$$||x^{*}||=||x||.$$
This map $||\cdot||:LS(M)\rightarrow L^{0}(\Omega)$ is called the center-valued norm on $LS(M).$

Finally we recall the following main result of the paper  \cite{Alb2}:

\textbf{Theorem 2.1.} \emph{Let  $M$ be a type $I$ von Neumann algebra with the center $Z$ and
a faithful normal semi-finite trace
$\tau.$ Then every $Z$-linear derivation $D$ on the algebra $S_0(M,
\tau)$ is spatial, namely
$$D(x)=ax-xa,\,x\in S_0(M, \tau)$$ for an appropriate operator $a\in S(M, \tau).$}

\begin{center} \textbf{3. Derivations on the algebra $S_0(M,
\tau)$ }\end{center}

In this section we give the main result of the paper and describe derivations on the algebra
$S_0(M, \tau)$ of all $\tau$-compact operators for type I von
Neumann algebra $M$ with a faithful normal semi-finite trace $\tau.$
 We shall do it step by step, considering separate cases.

\textbf{A. The case of finite von Neumann algebras of type  I. }

Let $N$ be a commutative von Neumann algebra, then $N\cong
L^{\infty}(\Omega)$ for an appropriate measure space
$(\Omega,\Sigma,\mu).$ It has been proved in  \cite{Ber},
\cite{Kus1} that the algebra $LS(N)=S(N)\cong L^{0}(\Omega)$
admits non trivial derivations if and only if the measure space
$(\Omega,\Sigma,\mu)$ is not atomic.

Let  $\tau$ be a faithful normal semi-finite trace on the commutative von
Neumann algebra $N$ and suppose that the Boolean
algebra $P(N)$ of projections is not atomic. This means that there exists a projection  $z\in N$ with
$\tau(z)<\infty$ such that the Boolean algebra of projection in $zN$ is continuous (i.e. has no atom). Since
$zS_0(N, \tau)=zS(N)\cong
S(zN),$ the algebra $zS_0(N, \tau)$ admits a non trivial derivation  $D$ (see \cite{Ber}, \cite{Kus1}). Putting
$$D_0(x)=D(zx),\, x\in S_0(N, \tau)$$
we obtain a non trivial derivation on the algebra $S_0(N, \tau).$ Therefore,
we have that if a commutative von Neumann
algebra $N$
has a non atomic Boolean algebra of projections then the algebra  $S_0(N, \tau)$ admits a non zero derivation.

Given an arbitrary derivation $D:S_{0}(N,\tau) \rightarrow S_{0}(N,\tau)$ the element
$$z_D=\inf\{z\in P(N): zD=D\}$$
is called the support of the operator  $D.$

\textbf{Lemma 3.1.} \emph{If $N$ is a commutative von Neumann algebra with a
faithful normal semi-finite trace  $\tau$
and $D:S_{0}(N,\tau) \rightarrow S_{0}(N,\tau)$ is a derivation, then $\tau(z_{D})<\infty.$}

Proof. Suppose the opposite, i.e.  $\tau(z_{D})=\infty.$ Then
there exists a sequence of mutually orthogonal projections
$z_{n}\in N,\,n=1,2...,$ with $z_{n}\leq z_{D},\ 1\leq
\tau(z_{n})<\infty.$ For $z=\sup\limits_{n}{z_{n}}$ we have
$\tau(z)=\infty.$ Since $\tau(z_{n})<\infty$  for all $n=1,2...,$
it follows that $z_{n}S_{0}(N,\tau)=z_n S(N)\cong S(z_{n}N).$
Define a derivation $D_{n}: z_{n}S(N)\rightarrow z_{n}S(N)$ by
$$D_{n}(x)=z_{n}D(x), \ x\in z_{n}S(N).$$
Since $z_{n}S(N)\cong S(z_{n}N)$ and $z_{D_{n}}=z_{n}$,
\cite[Lemma 3.2]{Alb4} implies that for each $n\in\mathbb{N}$
there exists an element $\lambda_{n}\in z_{n}N$ such that
$|\lambda_{n}|\leq n^{-1}z_{n}$ and $|D_{n}(\lambda_{n})|\geq
z_{n}.$

Put $\lambda=\sum \lambda_{n}. $
Then $|\lambda|\leq\sum\limits_{n}n^{-1}z_n$ and therefore $\lambda\in S_{0}(N,
\tau).$ On other hand
$$|D(\lambda)|=|D(\sum \lambda_{n})|=|D(\sum z_n\lambda_{n})|=
|\sum z_n D(\lambda_{n})|=\sum |D_n(\lambda_{n})|\geq \sum z_n=z,$$
i.e. $|D(\lambda)|\geq z.$ But $\tau(z)=\infty$, i.e. $z\notin S_{0}(N, \tau).$ Therefore
$D(\lambda)\notin S_{0}(N, \tau).$ The contradiction shows that $\tau(z_{D})<\infty.$ The proof is complete. $\blacksquare$

For general (non commutative) finite type I von Neumann algebras
we have

 \textbf{Lemma 3.2.} \emph{Let $M$ be a finite von Neumann
algebra of type I with the center  $Z$ and let
$D:S_{0}(M,\tau)\rightarrow S_{0}(M,\tau)$ be a derivation. If
$D(\lambda)=0$ for every $\lambda$ from the center
$Z(S_{0}(M,\tau))$ of $S_{0}(M,\tau),$ then $D$ is $Z$-linear.}

Proof. Take  $\lambda\in Z$ and choose a central projection $z$ in
$M$ with  $\tau(z)<\infty.$ Since  $z,\,z\lambda\in
Z(S_{0}(M,\tau)),$ we have that  $D(z)=D(z\lambda)=0.$

For  $x\in S_{0}(M,\tau)$ one has  $$D(z\lambda x)=D(z\lambda)x+z
\lambda D(x)=z\lambda D(x),$$ i.e.
$$D(z\lambda
x)=z\lambda D(x).$$ On the other hand
$$D(z\lambda x)=D(z)\lambda x+zD(\lambda x)=zD(\lambda x),$$
i.e.
$$D(z\lambda
x)=zD(\lambda x).$$ Therefore $zD(\lambda x)=z\lambda D( x). $
Since  $z$ is an arbitrary with  $\tau(z)<\infty$ this implies
(taking $z\uparrow \textbf{1}$) that $D(\lambda x)=\lambda D( x)$
for all $\lambda\in Z$ and $x\in S_0(M, \tau),$ i.e. $D$ is
$Z$-linear. The proof is complete. $\blacksquare$

Now let  $M$ be a homogeneous von Neumann algebra of type $I_{n},
n \in \mathbb{N}$,  with the center $Z$ and with a faithful normal
semi-finite trace  $\tau.$ Then  the algebra $M$ is *-isomorphic
with the algebra $M_n(Z)$ of all  $n\times n$- matrices over $Z,$
and the algebra  $S_0(M, \tau)$ is *-isomorphic with the algebra
 $M_n(S_0(Z, \tau_Z))$ of all  $n\times n$ matrices over
$S_0(Z, \tau_Z),$ where  $\tau_Z$ is the restriction of the trace
 $\tau$ onto the center  $Z.$

 If  $e_{i,j},\,i,j=\overline{1, n},$
are the matrix units in  $M_n(S_0(Z, \tau_Z)),$ then each element
$x\in M_n(S_0(Z, \tau_Z))$ has the form
 $$x=\sum\limits_{i,j=1}^{n}\lambda_{i,j}e_{i,j},\,\lambda_{i,j}\in S_0(Z, \tau_Z),\,i,j=\overline{1, n}.$$

Let  $\delta:S_0(Z, \tau_Z)\rightarrow S_0(Z, \tau_Z)$ be a
derivation and let  $z_\delta$ be its support. By Lemma 3.1 we
have  $\tau(z_\delta)<\infty,$ and therefore $z_\delta S_0(M,
\tau)=z_\delta S(M).$ Thus given an arbitrary
$\sum\limits_{i,j=1}^{n}\lambda_{i,j}e_{i,j}\in M_n(S_0(Z,
\tau_Z))$ one has
$$\sum\limits_{i,j=1}^{n}\delta(\lambda_{i,j})e_{i,j}=
z_\delta\sum\limits_{i,j=1}^{n}\delta(\lambda_{i,j})e_{i,j}\in
z_\delta M_n(S(Z))=z_\delta M_n(S_0(Z, \tau_Z))\subset M_n(S_0(Z,
\tau_Z)).$$ This implies that by putting
 $$D_{\delta}(\sum\limits_{i,j=1}^{n}\lambda_{i,j}e_{i,j})=
 \sum\limits_{i,j=1}^{n}\delta(\lambda_{i,j})e_{i,j} \eqno (1)$$
 we obtain a well-defined linear operator
 $D_\delta$ on the algebra $M_n(S_0(Z, \tau_Z)).$ Moreover
 $D_\delta$ is a derivation on the algebra  $M_n(S_0(Z, \tau_Z))$
 and its restriction onto the center of the algebra  $M_n(S_0(Z,
 \tau_Z))$ coincides with the given $\delta.$

 Now let the trace $\tau$ be finite. Then  $S_0(M, \tau)=S(M, \tau)=S(M).$ Consider a family
    $\{e_i\}_{i=1}^{n}$ of mutually orthogonal and mutually equivalent abelian projections in
    the von Neumann algebra  $M$ of type  $I_n,\, n\in\mathbb{N}.$
Put  $e=\sum\limits_{i=1}^{n-1}e_i.$ Then  $eMe$ is a von Neumann
algebra of type  $I_{n-1},$ and
$$S_0(Z, \tau_Z)\cong Z(eS_0(M, \tau)e)\cong Z(S_0(M, \tau)).$$

\textbf{Remark 1.} From now on we shall identify these isomorphic
abelian von Neumann algebras. In this case the element $\lambda$ from $S_0(Z, \tau_Z)$  corresponds to $\lambda e$
from $Z(eS_0(M, \tau)e)$ and to $\lambda\textbf{1}$ from
$Z(S_0(M, \tau)).$

Consider a derivation $D$ on the algebra  $S_0(M, \tau).$ Since
$D$ maps $Z(S_0(M, \tau))$ into itself, its
restriction  $D|_{Z(S_0(M, \tau))}$ induces a derivation
 $\delta$ on $S_0(Z, \tau_Z)\cong  Z(S_0(M, \tau)),$ i.e.
 $$D(\lambda\textbf{1})=\delta(\lambda)\textbf{1},\,\lambda\in S_0(Z, \tau_Z).$$
 Let  $D_e$ be the derivation on $eS_0(M, \tau)e$ defined as
 $$D_e(x)=eD(x)e, \,x\in eS_0(M, \tau)e.$$
 Since  $Z(eS_0(M, \tau)e)\cong Z(S_0(M, \tau)),$ the restriction of
 $D_e$ onto  $Z(eS_0(M, \tau)e)$ also generates a derivation, denoted by
  $\delta_e,$ on $S_0(Z, \tau_Z),$ i.e.
 $$D_e(\lambda e)=\delta_e(\lambda)e,\,\lambda\in S_0(Z, \tau_Z).$$

\textbf{Lemma 3.3.} \emph{The derivations $\delta$ and $\delta_e$ on
$S_0(Z, \tau_Z)$ coincide.}

 Proof.  Since  $e$ is a projection it is clear that  $eD(e)e=0$ and therefore
  $$\delta_e(\lambda)e=D_e(\lambda e)=eD(\lambda e)e=
 eD(\lambda\textbf{1})e+e\lambda D(e)e=eD(\lambda\textbf{1})e=\delta(\lambda)e,$$ i.e.
 $$\delta_e(\lambda)e=\delta(\lambda)e$$
 for any  $\lambda\in S_0(M, \tau).$ Therefore (see Remark 1)
 $\delta_e(\lambda)=\delta(\lambda),$ i.e.
  $\delta\equiv \delta_e.$ The proof is complete. $\blacksquare$

The following lemma describes derivations on the algebra of
$\tau$-compact operators for type  $I_n§ §(n\in\mathbb{N})$ von
Neumann algebras.

 \textbf{Lemma 3.4.} \emph{Let  $M$ be a homogenous von Neumann algebra of type
  $I_{n}, n \in \mathbb{N}$, with a faithful normal semi-finite
  trace $\tau.$ Every derivation  $D$ on the algebra $S_0(M, \tau)$ can be uniquely represented as a sum
  $$D=D_{a}+D_{\delta ,}$$ where  $D_{a}$ is a spatial derivation implemented by an element  $a\in S(M, \tau)$
while $D_{\delta} $ is the derivation of the form (1) generated by
a derivation $\delta$ on the center of $S_0(M, \tau)$ identified
with $S_0(Z, \tau_Z)$.}

Proof. Let  $D$ be an arbitrary derivation on the algebra $S_0(M,
\tau)\cong M_n(S_0(Z, \tau_Z)).$ Consider its restriction $\delta$
onto the center  $S_0(Z, \tau_Z)$ of this algebra, and let
$D_\delta$ be the derivation on the algebra $M_n(S_0(Z, \tau_Z))$
constructed as in (1). Put $D_1=D-D_\delta.$ Given any  $\lambda\in
S_0(Z, \tau_Z)$ we have
$$D_1(\lambda)=D(\lambda)-D_\delta(\lambda)=D(\lambda)-D(\lambda)=0,$$
i.e. $D_1$ is identically zero on  $S_0(Z, \tau_{Z}).$ By Lemma
3.2 $D_1$ is $Z$-linear and by theorem 2.1 we obtain that  $D_1$
is  spatial derivation and thus $D_1=D_a$ for an appropriate $a\in
M_n(S(Z, \tau_Z)).$ Therefore $D=D_a+D_\delta.$

Suppose that
$$D=D_{a_{1}}+D_{\delta_{1}}=D_{a_{2}}+D_{\delta_{2}}.$$ Then
$D_{a_{1}}-D_{a_{2}}=D_{\delta_{2}}-D_{\delta_{1}}.$ Since
$D_{a_{1}}-D_{a_{2}}$ is identically zero on the center of the
algebra $M_n(S_0(Z, \tau_Z))$ this implies that
$D_{\delta_{2}}-D_{\delta_{1}}$ is also identically zero on the
center of  $M_n(S_0(Z, \tau_Z)).$ This means that
$\delta_{1}=\delta_{2},$ and therefore  $D_{a_{1}}=D_{a_{2}},$
i.e. the decomposition of $D$ is unique. The proof is complete. $\blacksquare$

Now let  $M$ be an arbitrary finite von Neumann algebra of type I
with the center $Z.$ There exists a family  $\{z_n\}_{n\in F},$
$F\subseteq\mathbb{N},$ of central projections from $M$ with
$\sup\limits_{n\in F}z_n=\textbf{1}$ such that the algebra  $M$ is
*-isomorphic with the  $C^{*}$-product of von Neumann algebras
$z_n M$ of type  I$_{n}$ respectively, $n\in F,$ i.e.
$$M\cong\bigoplus\limits_{n\in F}z_n M.$$
In this  case we have that
$$S_0(M, \tau)\subseteq\prod\limits_{n\in F}S_0(z_n M, \tau_n),$$
where  $\tau_n$ is the restriction of the trace $\tau$ onto $z_n
M,\,n\in F.$

Suppose that    $D$ is a derivation on  $S_0(M, \tau),$ and
$\delta$ is its restriction onto its center  $S_0(Z, \tau_Z).$
Since  $\delta$ maps each  $z_nS_0(Z, \tau_Z))\cong Z(S_0(z_n M,
\tau_n))$ into itself,  $\delta$ generates a derivation $\delta_n$
on $z_nS_0(Z, \tau_Z)$ for each $n\in F.$

Let     $D_{\delta_n}$ be the derivation on the matrix algebra
$M_n(z_nZ(S_0(M, \tau)))\cong S_0(z_nM, \tau_n)$ defined as in
(1). Put
$$D_\delta(\{x_n\}_{n\in F})=\{D_{\delta_n}(x_n)\},\,\{x_n\}_{n\in F}\in S_0(M, \tau).\eqno(2)$$
Then the map  $D$  is a derivation on $S_0(M, \tau).$

Now Lemma 3.4 implies the following result:

\textbf{Lemma 3.5.} \emph{Let  $M$ be a finite von Neumann algebra
of type I with a faithful normal semi-finite trace $\tau.$ Each
derivation  $D$ on the algebra  $S_0(M, \tau)$ can be uniquely
represented in the form
$$D=D_{a}+D_{\delta ,}$$
where  $D_{a}$ is a spatial derivation implemented by an element
$a\in S(M, \tau),$ and $D_{\delta} $ is a derivation given as
(2).}

\newpage

\textbf{B. The case of type  I$_{\infty}$ von Neumann algebras. }

We shall start the consideration of this case by the description
of the center of the algebra of $\tau$-compact operators with
respect to type  I$_{\infty}$ von Neumann algebras.

\textbf{Lemma 3.6.} \emph{Let  $M$ be a type  $I_{\infty}$ von
Neumann algebra with the center $Z.$ Then}

    \emph{a) the centers of the algebras $S(M)$ and  $S(M,\tau)$ coincide with $Z;$}

    \emph{b) the center of the algebra  $S_{0}(M,\tau)$ is trivial, i.e.  $Z(S_{0}(M,\tau))=\{0\}.$}

   Proof. a) Suppose that  $z\in S(M), \ z\geq 0,$  is a central element and let
     $z=\int\limits_{0}^{\infty}\lambda\,de_{\lambda}$ be its spectral resolution. Then  $e_{\lambda}\in Z$ for all
     $\lambda>0.$ Assume that  $e_{n}^{\bot}\neq 0$ for all  $n\in\mathbb{N}.$
     Since $M$ is of type  I$_{\infty},$ $Z$ does not contain non-zero finite projections. Thus
   $e_{n}^{\bot}$ is infinite for all  $n\in\mathbb{N},$ which
   contradicts the condition  $z\in S(M).$ Therefore there
   exists  $n_{0}\in\mathbb{N}$ such that  $e_{n}^{\bot}=0$ for all  $n\geq n_0,$
  i.e. $z\leq n_0\textbf{1}.$ This means that   $z\in Z,$ i.e. $Z(S(M))=Z.$ Similarly $Z(S(M, \tau))=Z.$

  b) Let $z\in Z(S_0(M, \tau)), \ z\geq 0.$ Take a projection $p\in
  M$ with $\tau(p)<\infty.$ Then $p\in S_0(M, \tau)$ and therefore  $zp=pz.$
Since  $M$ is semi-finite this implies that $zp=pz$ for all $p\in
P(M).$ Since the linear span of  $P(M)$ is dense in $S(M, \tau)$
in the measure topology, we have that  $zx=xz$ for all $x\in S(M,
\tau),$ i.e. $z\in Z(S(M, \tau))=Z.$

    Suppose that
    $z=\int\limits_{0}^{\infty}\lambda\,de_{\lambda}$ is the
    spectral resolution of $z.$ Then
     $e_{\lambda}\in Z$ for all  $\lambda>0.$ Since
$z\in S_0(M, \tau)$ we have that $e_{\lambda_0}^{\perp}$ is a
finite projection for all $\lambda_0>0.$ But $M$ does not contain
any non zero central finite projection, because it is of type
 I$_{\infty}.$ Therefore  $e_{\lambda}^{\perp}=0$ for all
 $\lambda>0,$ i.e.  $z=0.$ Thus  $Z(S_0(M,
    \tau))=\{0\}.$ The proof is complete. $\blacksquare$

It should be noted that the center of the algebra $LS(M),$ for a
general von Neumann algebra  $M$  coincides with $LS(Z)$ and thus
contains $Z.$ This was an essential point in the proof of theorems
concerning the description of derivations on the algebra $LS(M)$
of locally measurable operators with respect to type I von Neumann
algebras $M$ (see [4]). Lemma 3.6 shows that this is not the case
for the algebra  $S_{0}(M,\tau)$ because the center of this
algebra may be trivial. Thus the methods of the paper \cite{Alb4}
can not be directly applied for the description of derivations of
the algebra of  $\tau$-compact operators with respect to type I
von Neumann algebras.

We are now in position to prove one of the main results of the
paper.

 \textbf{Theorem 3.7.} \emph{If $M$ is a type  $I_\infty$ von Neumann algebra with a faithful normal semi-finite trace
  $\tau,$ then every derivation on the algebra $S_0(M, \tau)$ is spatial and implemented by an element of the algebra
   $S(M, \tau).$ }

 The proof of the theorem consists of several lemmata.

\textbf{Lemma  3.8.} \emph{Let   $z\in Z$ be a central projection
from $M$ and let  $x\in S_0(M, \tau).$ Then}
 $$D(zx)=zD(x).$$

Proof. Without loss of generality we may suppose that  $x\geq
0,$ i.e. $x=y^{2}$ for some $y\in S_0(M, \tau).$ From the Leibniz
rule for derivations we obtain
$$D(zx)=D(zyzy)=D(zy)zy+zyD(zy)=
z[D(zy)y+yD(zy)].$$ Therefore
$$z^{\perp}D(zx)=0.$$
Similarly we have that
$$zD(z^{\perp}x)=0.$$
Further
$$zD(x)=zD((z+z^{\perp})x)=zD(zx)+zD(z^{\perp}x)=zD(zx),$$
i.e.
$$zD(x)=zD(zx).$$
On the other hand
$$D(zx)=(z+z^{\perp})D(zx)=zD(zx)+z^{\perp}D(zx)=zD(zx),$$
i.e.
$$D(zx)=zD(zx).$$
 Therefore   $$D(zx)=zD(x).$$ The proof is complete.

  \textbf{Lemma  3.9.} \emph{Suppose that $\lambda\in Z,$
  $p\in P(M),\, \tau(p)<\infty.$ Put $y=D(\lambda p)-\lambda D(p).$ Then}
 $$p^{\perp}yp^{\perp}=0.$$

 Proof. From $$D(p)=D(pp)=D(p)p+pD(p)$$ and
  $$D(\lambda p)=D(\lambda pp)=D(\lambda p)p+\lambda pD(p)$$ we
  obtain
 $$p^{\perp}D(\lambda p)p^{\perp}=p^{\perp}\lambda D(p)p^{\perp}=0$$ and in particular $p^{\perp}yp^{\perp}=0.$
The proof is complete. $\blacksquare$

 \textbf{Lemma  3.10.} \emph{For each  $\lambda\in Z$ and for every abelian projection
  $p\in P(M)$ with $\tau(p)<\infty$ we have}
 $$D(\lambda p)=\lambda D(p).$$

Proof.  Let $z$ be the central cover of the projection $p,$ (i.e.
the least central projection majorating $p$). Lemma 3.8 implies
that the derivation $D$ maps the algebra $zS_0(M, \tau)$ into
itself. Therefore passing if necessary to the algebra $zM$ and to
the derivation  $zD$ we may assume without loss of generality
that $z=\textbf{1},$ i.e. that $p$ is a faithful projection. Take
an arbitrary projection  $p_0$ with the central cover $\textbf{1}$
such that $p_0\leq p^{\perp}$ and such that the von Neumann
algebra $p_0Mp_0$ is of type  I$_{\aleph_{0}},$ where $\aleph_{0}$
is the countable cardinal number. Then there exists a sequence of
mutually orthogonal and pairwise equivalent abelian projections
$\{p_{n}\}_{n=2}^{\infty}$ in $M$ with $\sum\limits_{n=2}^{\infty}
p_{n}=p_0.$ Putting $p_1=p$ we obtain that projections $p_1$ and $p_2$ are equivalent ($p_1\sim p_n$)  and thus
 $\tau(p_n)=\tau(p_1)<\infty$ for all  $n\geq 2.$

Put $e_{n}=\sum\limits_{k=1}^{n} p_{k},\ n\geq 1$. Then
$e_{n}Me_{n}$ is a homogeneous von Neumann algebra of type
$I_{n},\ n\geq 1,$ and the restriction $\tau_{n}$ of the trace
$\tau$ onto $e_{n}Me_{n}$ is finite, and therefore  $e_{n}S_0(M,
\tau)e_{n}\cong S(e_{n}Me_{n}),\,n=1,2....$

Define a derivation $D_n$ on $e_{n}S_0(M, \tau)e_{n}$ as
follows
$$D_{n}(x)=e_{n}D(x)e_{n},\, x \in e_{n}S_0(M,
\tau)e_{n}.$$ By lemma  3.4 for each  $n$ there exist an element
$a_{n} \in e_{n}S_0(M, \tau)e_{n}$ and a derivation $\delta_{n}$
on  $e_{1}S_0(M, \tau)e_{1}$ identified with $Z(e_{n}S_0(M,
\tau)e_{n})$ (see Remark 1) such that
$$D_{n}=D_{a_{n}}+D_{\delta_{n}}.\eqno (3)$$
Since  $D_n=e_n D_{n+1}e_n$ lemma 3.3 implies that
$\delta_{n}=\delta_{n+1},\ n\geq 1.$ Denote $\delta=\delta_n.$

Given a sequence  $\Lambda=\{\lambda_{n}\}$ in  $Z$ with
$|\lambda_{n}|\leq \frac{\textstyle 1}{\textstyle n}\textbf{1},\ n
\in \mathbb{N},$ put
$$x_{\Lambda}=\sum\limits_{n=1}^{\infty}\lambda_{n}p_{n}.$$
Let us show that  $x_{\Lambda}\in S_0(M, \tau).$

For an arbitrary  $\varepsilon>0$ there exists  $n_{0}\in
\mathbb{N}$ such that $\frac{\textstyle 1}{\textstyle
n_0}<\varepsilon.$ Put
$$p_\varepsilon=\textbf{1}-\sum
\limits_{n=1}^{n_{0}-1}p_{n},$$ then
 $\tau(p^{\perp}_\varepsilon)=\tau(\sum
\limits_{n=1}^{n_{0}-1}p_{n})=(n_{0}-1)\tau(p_1)<\infty.$ Moreover
$$||x_{\Lambda}p_\varepsilon||_{M}=||\sum
\limits_{n=n_{0}}^{\infty}\lambda_{n}p_{n}||_{M}=\sup\limits_{n\geq
n_0}||\lambda_{n}||_{M}\leq {1 \over n_{0}}<\varepsilon.$$
 This means that  $x_{\Lambda}\in
S_{0}(M,\tau).$ For each  $n \in \mathbb{N}$ we have
$$x_{\Lambda}p_{n}=p_{n}x_{\Lambda}=\lambda_{n} p_{n}.$$
By the Leibniz identity for derivations this implies
$$D(\lambda_{n}p_{n})=D(x_{\Lambda})p_{n}+x_{\Lambda}D(p_{n}).$$
Multiplying by $p_{n}$ from  both sides we obtain
$$p_{n}D(\lambda_{n}p_{n})p_{n}=p_{n}D(x_{\Lambda})p_{n}p_{n}+p_{n}x_{\Lambda}D(p_{n})p_{n},$$
i.e.
$$p_{n}D(\lambda_{n}p_{n})p_{n}=p_{n}D(x_{\Lambda})p_{n}+\lambda_{n}p_{n}D(p_{n})p_{n},$$
Therefore
$$p_{n}D(\lambda_{n}p_{n})p_{n}=p_{n}D(x_{\Lambda})p_{n}\eqno
(4)$$ because $pD(p)p=0$ for every projection $p\in M.$

On the other hand
$$p_{n}D(\lambda_{n}p_{n})p_{n}=p_{n}e_{n}D(\lambda_{n}p_{n})e_{n}p_{n}=p_{n}D_{n}(\lambda_{n}p_{n})p_{n}.$$
From  (3) we obtain
$$p_{n}D(\lambda_{n}p_{n})p_{n}=p_{n}D_{a_{n}}(\lambda_{n}p_{n})p_{n}+p_{n}D_{\delta}(\lambda_{n}p_{n})p_{n}.$$
Since  $D_{a_{n}}$ is a spatial derivation (and hence it is
$Z$-linear), we have that
$$p_{n}D_{a_{n}}(\lambda_{n}p_{n})p_{n}=\lambda_{n}p_{n}D_{a_{n}}(p_{n})p_{n}=0.$$
From
$$p_{n}D_{\delta}(\lambda_{n}p_{n})p_{n}=\delta(\lambda_{n})p_{n},$$
we obtain
$$p_{n}D(\lambda_{n}p_{n})p_{n}=\delta(\lambda_{n})p_{n}.\eqno
(5)$$ Now  (4) and  (5) imply
$$p_{n}D(x_{\Lambda})p_{n}=\delta(\lambda_{n})p_{n}.$$

Suppose that  $\delta\neq 0.$ Then  \cite[Lemma 3.2]{Alb4} implies
the existence of a sequence  $\Lambda=\{\lambda_{n}\}$ in $Z$ with
 $|\lambda_{n}|\leq \frac{\textstyle 1}{\textstyle n}\textbf{1},\
n \in \mathbb{N},$ and a projection  $\pi\in Z,\,\pi\neq0$ such
that
$$|\delta(\lambda_{n})|\geq n\pi, \ n\in \mathbb{N}.$$

Further
$$||p_{n}D(x_{\Lambda})p_{n}||\leq||p_{n}|| ||D(x_{\Lambda})|| ||p_{n}||=||D(x_{\Lambda})||$$
and $$||\delta(\lambda_{n})p_{n}||=|\delta(\lambda_{n})|.$$
Therefore  $$||D(x_{\Lambda})|| \geq |\delta(\lambda_{n})|,$$ and
hence $$||D(x_{\Lambda})|| \geq \pi n, \ n\geq 2.$$

The  last inequality contradicts the choice of $\pi\neq 0.$
Therefore  $\delta\equiv 0,$ i.e. from (3) we obtain that
$D_{n}=D_{a_{n}}.$ Since  $D_{a_{n}}$ is a spatial derivation and
the center of the algebra  $Z(e_nMe_n)$ coincides with $e_nZ,$ it
follows that $D_n$ is  $e_nZ$-linear. Thus
$$D_{n}(\lambda e_npe_n)=\lambda e_nD_{n}(e_npe_n)\eqno (6)$$
for all  $\lambda\in Z.$ Since the projection  $e_n$ is in $S_0(M,
\tau)$ and it commutes with $p$ we have  $$D_n(e_n pe_n)=D_n(e_n
p)=e_nD(e_n p)e_n=e_nD(e_n)pe_n+e_nD(p)e_n=$$$$=e_nD(e_n)e_n
p+e_nD(p)e_n=e_nD(p)e_n,$$ i.e.
$$\lambda D_n(e_n pe_n)=\lambda e_nD(p)e_n.\eqno (7)$$
In a similar way we obtain
 $$D_n(\lambda e_n pe_n)=e_nD(\lambda p)e_n.\eqno (8)$$
Now (6), (7) and  (8) imply
$$e_{n}D(\lambda  p)e_{n}=e_{n}\lambda D(p)e_{n}$$ for all $n\in
\mathbb{N}$.

Set  $y=D(\lambda p)-\lambda D(p).$ Then  $e_n ye_n=0.$ From
$e_1=p_1=p,$  we have   $pyp=0.$ By Lemma 3.9 we have
$p^{\perp}yp^{\perp}=0.$ Multiplying the equality  $e_n ye_n=0$ by
 $p$ from the left side we obtain  $pye_n=0$ for all $n\in\mathbb{N}.$
 Since   $e_n\uparrow p_0+p,$ it follows that
 $py(p_0+p)=0,$ i.e.  $pyp_0=0.$ Since  $p_0$ is an
 arbitrary projection with the central cover $\textbf{1}$
such that $p_0\leq p^{\perp}$ and such that the von Neumann
algebra $p_0Mp_0$ is of type  I$_{\aleph_{0}},$  we obtain that
$pyp^{\perp}=0.$

 Similarly  $p^{\perp}yp=0.$ Therefore
 $$pyp=pyp^{\perp}=p^{\perp}yp=p^{\perp}yp^{\perp}=0$$ and hence
 $$y=pyp+pyp^{\perp}+p^{\perp}yp+p^{\perp}yp^{\perp},$$ i.e. $D(\lambda
p)=\lambda D(p).$ The proof is complete.

\textbf{Lemma  3.11.} \emph{Suppose that   $\lambda\in Z$ and
$x\in S_0(M, \tau).$ Then  }
 $$D(\lambda x)=\lambda D(x).$$

Proof. Case (i).  $x=p$ is a projection and
$$p=\sum\limits_{i=1}^{k}p_i,\eqno (9) $$ where  $p_i,\, i=\overline{1, k}$
are mutually orthogonal abelian projections with
$\tau(p_i)<\infty.$ By Lemma  3.10 we have  $D(\lambda
p_i)=\lambda D(p_i).$ Therefore
$$D(\lambda p)=D(\lambda\sum\limits_{i=1}^{k}p_i)=\sum\limits_{i=1}^{k}D(\lambda
p_i)=\sum\limits_{i=1}^{k}\lambda D(p_i)=\lambda
D(\sum\limits_{i=1}^{k}p_i)=\lambda D(p),$$ i.e.
$$D(\lambda p)=\lambda D(p).$$

Case (ii).  $x=p$ is a projection with $\tau(p)<\infty.$ Then
$pMp$ is a finite von Neumann algebra of type I, and therefore
there exists a sequence of mutually orthogonal central projections
 $\{z_n\}$ such that  each $p_n=z_n p$ is a projection of the form (9). From the above case we have
 $D(\lambda p_n)=\lambda D(p_n).$ This and Lemma  3.8 imply that
 $$z_n D(\lambda p)=D(\lambda z_n p)=D(\lambda p_n)=\lambda D(p_n)=\lambda D(z_n p)=\lambda z_n D(p).$$
i.e.
 $$z_n D(\lambda p)=z_n\lambda  D(p)$$ for all $n.$
Therefore
$$D(\lambda p)=\lambda  D(p)$$

Case (iii). Let  $x\in S_0(M, \tau)$ be an element such that
$xp=x$ for some projection $p$ with $\tau(p)<\infty.$ Then
$$D(\lambda x)=D(\lambda x p)=D(x \lambda p)=D(x)\lambda p+xD(\lambda
p)=$$$$=D(x)\lambda p+x\lambda D(p)=\lambda(D(x)p+xD(p)=\lambda
D(xp)=\lambda D(x),$$ i.e. $D(\lambda x)=\lambda D(x).$

Case (iv).   $x$ is an arbitrary element from $S_0(M, \tau).$
Take a projection $p$ with  $\tau(p)<\infty.$ Put $x_0=xp.$ From
the case (iii) we have  $D(\lambda x_0)=\lambda D(x_0).$ Now one
has
$$D(\lambda x_0)=D(\lambda x p)=D(\lambda x)p+\lambda xD(p),$$
i.e.
$$D(\lambda x)p=D(\lambda x_0)-\lambda xD(p).$$
On the other hand
$$D(\lambda x_0)=\lambda D(x_0)=\lambda D(xp)=\lambda D(x)p+\lambda x D(p),$$
i.e.
$$\lambda D(x)p=D(\lambda x_0)-\lambda x D(p).$$
Therefore   $\lambda D(x)p=D(\lambda x)p.$ Since  $p$ is an
arbitrary with $\tau(p)<\infty,$ this implies
$$D(\lambda x)=\lambda D(x).$$ The proof is complete. $\blacksquare$

Proof of Theorem 3.7.

By Lemma  3.11 the derivation  $D: S_0(M, \tau)\rightarrow S_0(M,
\tau)$ is  $Z$-linear. By Theorem  2.1 $D$ is spatial and moreover
$$D(x)=ax-xa,\,x\in S_0(M, \tau)$$ for an appropriate  $a\in S(M,
\tau).$ The proof is complete. $\blacksquare$

\textbf{C. The general case of type  I von Neumann algebras. }

Now we can obtain a complete description of derivations on the
algebra $S_0(M, \tau)$ of $\tau$-compact operators with respect to
a general type I von Neumann algebra $M$ with a faithful normal
semi-finite trace $\tau.$

Let $M$ be a type  I von Neumann algebra. There exists a central
projection  $z_0\in M$ such that

a) $z_0M$ is a finite von Neumann algebra;

b) $z_0^{\bot}M$ is a von Neumann algebra of type  I$_{\infty}.$

Now consider a derivation  $D$ on  $S_0(M, \tau)$ and let $\delta$
be its restriction onto its center  $Z(S_0(M, \tau)).$ By Lemma
3.6 we have $z_0^{\bot}Z(S_0(M, \tau))=\{0\},$ and therefore
$z_0^{\bot}\delta\equiv 0,$ i.e. $\delta=z_0\delta.$

Let   $D_\delta$ be the derivation on  $z_0S_0(M, \tau)$ defined
as in (2) and let us extend it to whole  on $S_0(M,
\tau)=z_0S_0(M, \tau)\oplus z_0^{\bot}S_0(M, \tau)$ by putting
$$D_\delta(x_1+x_2):=D_\delta(x_1),\,x_1\in z_0S_0(M, \tau),x_2\in z_0^{\bot}S_0(M,
\tau).\eqno (10)$$

The following theorem is the main result of the present paper, which gives
the general form of derivations on the algebra $S_0(M, \tau).$

\textbf{Theorem  3.12.} \emph{Let  $M$ be a type  $I$ von Neumann
algebra with a faithful normal semi-finite trace $\tau.$ Each
derivation  $D$ on  $S_0(M, \tau)$ can be uniquely represented in
the form
$$D=D_{a}+D_{\delta ,}\eqno(11)$$
where  $D_{a}$ is a spatial derivation implemented by an element
$a\in S(M, \tau),$ and $D_{\delta} $ is a derivation of the form
(10), generated by a derivation $\delta$ on the center of $S_0(M,
\tau)$.}

Proof. It  immediately follows from Lemma  3.5 and Theorem  3.7. $\blacksquare$

\textbf{D. An application to the description of cohomology
groups.}

Let  $A$ be an algebra. Denote by  $Der(A)$ the space of all
derivations (in fact it is a Lie algebra with respect to the
commutator), and denote  by $InDer(A)$ the subspace of all inner
derivations on $A$ (it is a Lie ideal in $Der(A)$).

 The factor-space   $H^{1}(A)=Der(A)/InDer(A)$ is called the first cohomology group of the algebra
  $A$ (see \cite{Dal}). It is clear that  $H^{1}(A)$ measures how much the space of all derivations
  on $A$ differs from the space on inner derivations.

Further we need the following property of the algebra of
$\tau$-compact operators from  \cite{Str}:
$$S(M, \tau)=M+S_0(M, \tau).\eqno (12)$$
Set  $C(M, \tau)=M\cap S_0(M, \tau)$ and consider  $M/(C(M,
\tau)+Z)$ -- the factor space of  $M$ with respect to the space
$C(M, \tau)+Z.$

For  $D_1, D_2\in Der(S_0(M, \tau))$ put
$$D_1 \sim D_2 \Leftrightarrow D_1-D_2\in InDer(S_0(M, \tau)).$$

Let  $D_1 \sim D_2.$ From Theorem 3.12 these derivations can be
represented in the form (11):
$$D_1=D_a+D_\delta,\, D_2=D_b+D_\sigma.$$
Since $D_1-D_2=D_c,$ where $c\in S_0(M, \tau)\subset S(M, \tau),$
from the uniqueness of a representation in the form (11) it follows
that  $D_a-D_b\in InDer(S_0(M, \tau))$ and $D_\delta=D_\sigma.$
Therefore  $\delta\equiv\sigma$ and
$$a-b\in
S_0(M, \tau)+Z(S(M, \tau)).\eqno (13)$$ According to  (12) we have
$$a=a_1+a_2, a_1\in M,\,a_2\in S_0(M, \tau),$$
$$b=b_1+b_2, b_1\in M,\,b_2\in S_0(M, \tau).$$
From  (13) it follows that
$$a_1-b_1\in (b_2-a_2)+S_0(M, \tau)+Z(S(M,
\tau))\subset S_0(M, \tau)+Z(S(M, \tau)).$$ Since  $a_1, b_1\in
M,$ we have that   $$a_1-b_1\in (S_0(M, \tau)+Z(S(M, \tau)))\cap
M\subset C(M, \tau)+Z$$ because $Z(S(M, \tau))\cap M=Z$ (cf. Lemma
3.6). Therefore
$$D_1 \sim D_2 \Leftrightarrow a_1-b_1\in C(M, \tau)+Z,\, \delta\equiv\sigma.$$

Thus we have following result.

\textbf{Theorem 3.13.} \emph{Let  $M$ be a type I von Neumann
algebra with the center $Z$ and with a faithful normal semi-finite
trace $\tau.$ Suppose that  $z_0$ is a central projection such
that  $z_0M$ is a finite von Neumann algebra, and  $z_0^{\bot}M$
is of type  I$_{\infty}.$ Then the group  $H^{1}(S_0(M, \tau))$ is
isomorphic with the group  $M/(C(M, \tau)+Z)\oplus H^{1}(S_0(z_0Z,
\tau_0)),$ where $\tau_0$ is the restriction of  $\tau$ onto
$z_0Z.$ In particular }

\emph{a) if the trace  $\tau$ is finite then $H^{1}(S_0(M,
\tau))\cong H^{1}(S_0(z_0Z, \tau_0));$}

\emph{b) if  $M$ is of type  I$_{\infty},$ then $H^{1}(S_0(M,
\tau))\cong M/(C(M, \tau)+Z).$ }

\vspace{1cm}

\textbf{Acknowledgments.} \emph{The second and third named authors
would like to acknowledge the hospitality of the $\,$ "Institut
f\"{u}r Angewandte Mathematik",$\,$ Universit\"{a}t Bonn
(Germany). This work is supported in part by the DFG 436 USB
113/10/0-1 project (Germany) and the Fundamental Research
Foundation of the Uzbekistan Academy of Sciences.}

\end{document}